\documentclass[11pt,oneside]{amsart}

\usepackage[T1]{fontenc}
\usepackage{amsmath,amssymb,amsthm,mathtools}
\usepackage{mathrsfs}
\usepackage{tikz-cd}
\usepackage{xcolor}
\usepackage{geometry}
\usepackage{hyperref}
\usepackage[nameinlink]{cleveref}

\DeclareMathAlphabet{\mathcal}{OMS}{cmsy}{m}{n}

\newcommand{\C}{\mathbb{C}}
\newcommand{\Q}{\mathbb{Q}}
\newcommand{\Z}{\mathbb{Z}}

\newcommand{\mc}{\mathcal}

\DeclareMathOperator{\Hom}{Hom}

\DeclareMathOperator{\Res}{Res}
\DeclareMathOperator{\ind}{ind}
\DeclareMathOperator{\Ind}{Ind}

\newcommand{\adm}{\mathrm{adm}}

\theoremstyle{plain}
\newtheorem{thm}{Theorem}[section]
\newtheorem*{mainthm}{Theorem A}
\newtheorem{lem}[thm]{Lemma}
\newtheorem{prop}[thm]{Proposition}

\theoremstyle{definition}
\newtheorem{defn}[thm]{Definition}

\theoremstyle{remark}
\newtheorem{rem}[thm]{Remark}

\numberwithin{equation}{section}

\title[Admissibility of unbounded complexes]
{Admissibility of unbounded complexes of smooth representations}

\author{Hargun Bhatia}
\address{Department of Mathematics, University of California San Diego,
9500 Gilman Drive, La Jolla, CA 92093-0112, USA}
\email{hpbhatia@ucsd.edu}
\date{November 11, 2026}
\hypersetup{
  pdftitle={Admissibility of unbounded complexes of smooth representations},
  pdfauthor={},
  colorlinks=true,
  linkcolor=blue,
  citecolor=blue,
  urlcolor=blue
}

\begin{document}

\begin{abstract}
Let $G$ be a $p$-adic Lie group and $k$ a field of characteristic $p$. We prove that a complex of smooth $k$-representations of $G$ has admissible cohomology if and only if it is globally admissible, i.e., its hypercohomology with respect to a torsion-free open compact pro-$p$ subgroup is finite dimensional in every degree. This extends the bounded-below and bounded-above comparison of Schneider--Sorensen to arbitrary unbounded complexes. 
\end{abstract}
\maketitle

\section{Introduction}\label{sec:introduction}
Many interesting smooth representations of $p$-adic Lie groups are infinite dimensional. A basic finiteness condition that is useful is that of admissibility. A smooth representation $V$ of a $p$-adic Lie group $G$ is called \textbf{admissible} if the space $V^K$ of vectors fixed by each compact open subgroup $K\subseteq G$ is finite dimensional. In the theory of representations over $\C$, any smooth irreducible representation of a $p$-adic reductive group is in fact admissible \cite[Théorème VI.2.2]{Renard}. Admissibility also underlies smooth duality: the canonical map
\[ V\longrightarrow  V^{\vee \vee}
\] is an isomorphism for an admissible smooth representation $V$ over $\C$, where $V^{\vee}$ denotes the smooth part of the linear dual of $V$ \cite[Proposition III.1.7]{Renard}.

In characteristic $p$, neither assertion holds in general. For $p>2$, Le constructed smooth irreducible representations of $\text{GL}_2(\Q_{p^3})$ which are not admissible, where $\Q_{p^3}$ denotes the degree $3$ unramified extension of $\Q_p$ \cite[Theorem 1.1]{Le2019}. Also, the smooth dual of an
infinite-dimensional irreducible admissible representation is zero
\cite[Proposition 3.9]{Kohlhaase2017}. Kohlhaase's higher smooth duals and their subsequent interpretation by Schneider--Sorensen in terms of derived smooth duals \cite{SchSor} provide a way to recover this information lost by just the ordinary smooth dual. 

Moreover, there is another natural reason to work with derived categories in characteristic $p$. Over $\C$, the functor of taking invariants $(-)^K$ under a compact open subgroup $K\subseteq G$ is exact, but it need not be exact in characteristic $p$. The cohomology functors $H^i(K,-)$ record this failure of exactness, and hence, passing to the derived category allows one to retain all this extra information. Thus, it's natural to ask how admissibility interacts with derived smooth duality. 

Let $G$ be a finite dimensional $p$-adic Lie group, $k$ be a field of characteristic $p$ and let $D(G)$ denote the unbounded derived category of smooth representations of $G$ in $k$-vector spaces. Schneider--Sorensen \cite{SchSor} construct the derived smooth duality functor
\[
\mathbb{D}=R\underline{\Hom}(-,k): D(G)^{\text{op}}\rightarrow D(G)
\]
and prove that the canonical biduality morphism
\[
\eta_{X^{\bullet}}: X^{\bullet}\longrightarrow \mathbb{D}^2(X^{\bullet})
\]
is an isomorphism for $X^{\bullet}\in D(G)$ whose cohomologies $h^i(X^{\bullet})$ are admissible $G$-representations for all $i\in \Z$. Such complexes are called \textit{admissible} and their full subcategory inside $D(G)$ is denoted by $D_{\text{adm}}(G)$. Schneider--Sorensen also define another notion of admissibility. Fix a torsion-free compact open pro-$p$ subgroup $U\subseteq G$. A complex $X^{\bullet}$ is \textit{globally admissible} if the hypercohomology $H^i(U, X^{\bullet})$ is finite dimensional for all $i\in \Z$. We denote its full subcategory inside $D(G)$ by $D(G)^a$. Schneider--Sorensen then prove that the canonical derived biduality morphism $\eta_{X^{\bullet}}$ is an isomorphism for $X^{\bullet}\in D(G)$ if and only if $X^{\bullet}\in D(G)^a$. This shows that the definition of globally admissible is independent of the choice of $U$ as well as that we have the inclusion 
\[
D_{\text{adm}}(G)\subseteq D(G)^a.
\]
Schneider--Sorensen further prove that the two subcategories above are equal for bounded-below or bounded-above complexes \cite[Proposition~4.9(i)]{SchSor}. Their argument for the bounded-below case proceeds inductively using standard truncations, beginning
with the lowest possibly nonzero cohomology degree. They also establish comparison for unbounded complexes with at least one vanishing differential \cite[Proposition~4.9(ii)]{SchSor}, and more generally for complexes that are split somewhere \cite[Remark~4.10]{SchSor}. 

The purpose of this paper is to establish the comparison for arbitrary unbounded complexes, with no splitting assumption. In fact, we prove a statement in each individual degree.
\begin{mainthm}
Let $U\subseteq G$ be a compact open pro-$p$ subgroup. Let $V^{\bullet}\in D(G)$ and $i\in \Z$. Then
\[
\dim_k H^i(U, V^{\bullet})< \infty \Rightarrow h^i(V^{\bullet})\; \text{is an admissible}\; G\text{-representation}.
\]
\end{mainthm}
Choosing $U$ torsion-free and combining this with the known inclusion then gives
\[
D_{\text{adm}}(G)=D(G)^a.
\]
The proof uses a homotopically minimal $K$-injective resolution of $\Res^G_U(V^{\bullet})$ in the sense of \cite[Appendix B]{Krause}. Minimality forces the differential to vanish on the $U$-invariants of such a resolution, which implies that $H^i(U, V^{\bullet})$ is identified with the $U$-invariants of the corresponding term of the resolution.

In Section \ref{sec:minimal_res}, we recall the necessary facts about homotopically minimal complexes and $K$-injective resolutions. In Section \ref{sec:proof}, we prove Theorem A and discuss its consequences. In the final section, we combine this result with Emerton’s admissibility theorem for higher ordinary parts to give another proof that the left adjoint of derived parabolic induction constructed by Heyer \cite{Heyer2023} preserves global admissibility.\\

\textbf{Acknowledgments:} I would like to thank my advisor, Claus Sorensen, for suggesting this problem and for his guidance throughout this work. I am grateful for his continued support and encouragement. I also thank Casey Perdue for listening carefully to my explanation of the argument and helping me check for errors. \\

\section{Homotopically minimal resolutions }\label{sec:minimal_res}
Let $\mc{A}$ be an abelian category. All complexes are cohomologically indexed, in particular, for a chain complex $X^{\bullet}$ in $\mc{A}$, $h^i(X^{\bullet}[m])=h^{i+m}(X^{\bullet})$ for all $i,m\in \Z$. We write $Z^n(X^{\bullet}):=\text{ker}(d^n:X^n\rightarrow X^{n+1})$ for all $n\in \Z$.

Recall that an inclusion $M\hookrightarrow N$ in $\mc{A}$ is called \textit{essential} if every nonzero subobject of $N$ has nonzero intersection with $M$. If furthermore, $N$ is an injective object, we say that $N$ is an \textit{injective envelope} of $M$.  We recall the notion of \textit{homotopically minimal} complexes as defined in \cite[Appendix B]{Krause}. 
\begin{defn}
A complex $X^{\bullet}$ in $\mc{A}$ is called \textbf{homotopically minimal} if any endomorphism $\phi:X^{\bullet}\rightarrow X^{\bullet}$ that is a chain homotopy equivalence is an isomorphism of chain complexes.    
\end{defn}
Clearly any complex with zero differentials is a homotopically minimal complex. There is a useful criterion for checking whether a complex is homotopically minimal. 

\begin{prop}[{\cite[Lemma B.1]{Krause}}]
Assuming that $\mc{A}$ admits injective envelopes, we have that a complex of injectives $X^{\bullet}$ is homotopically minimal if and only if the canonical map $Z^n(X^{\bullet})\hookrightarrow X^n$ is an injective envelope for all $n\in \Z$.  
\end{prop}

\textbf{Example.} If $\mc{A}=\text{Vect}_k$ is the category of vector spaces over a field $k$, then a complex $X^{\bullet}$ is homotopically minimal if and only if all the differentials are zero. Indeed, every vector space is injective and an inclusion of vector spaces is essential only when it's an equality. Hence, by the criterion above, we conclude that $Z^n(X^{\bullet})=X^n$ for all $n\in \Z\Rightarrow X^{\bullet}$ has zero differentials.

\textbf{Example.} Consider the complex of $k$-vector spaces: $[k\oplus k\xrightarrow[]{(x,y)\mapsto x}k]$. This is not homotopically minimal, however, we have a decomposition
\[
[k\oplus k\xrightarrow{(x,y)\mapsto x}k]
\cong
\underbrace{[k\longrightarrow 0]}_{\text{homotopically minimal}}
\oplus
\underbrace{[k\xrightarrow{\mathrm{id}}k]}_{\text{null-homotopic}}.
\]
Both two-term complexes are placed in degrees $0$ and $1$.

The above decomposition is more general. 
\begin{prop}[{\cite[Proposition B.2]{Krause}}]\label{krause_dec}
Let $\mc{A}$ be an abelian category admitting injective envelopes. Then every complex $X^{\bullet}$ with injective terms  has a decomposition 
$X^{\bullet}=X_1^{\bullet}\oplus X_2^{\bullet}$ such that $X_1^{\bullet}$ is homotopically minimal and $X_2^{\bullet}$ is nullhomotopic. Moreover, if there is another decomposition $X^{\bullet}=Y_1^{\bullet}\oplus Y_2^{\bullet}$ with $Y_1^{\bullet}$ homotopically minimal and $Y_2^{\bullet}$ nullhomotopic, the canonical map $X_1^{\bullet}\hookrightarrow X^{\bullet}\twoheadrightarrow Y_1^{\bullet}$ is an isomorphism of complexes. 
\end{prop}

Let $K(\mc{A})$ denote the homotopy category of chain complexes in $\mc{A}$. Recall that a complex $I^{\bullet}$ in $\mc{A}$ is called $K$-\textbf{injective} if $\Hom_{K(\mc{A})}(C^{\bullet}, I^{\bullet})=0$ for any acyclic complex $C^{\bullet}$. Clearly, being $K$-injective is preserved under direct summands. The full subcategory of complexes in $K(\mc{A})$ which are $K$-injective forms a triangulated subcategory denoted by $K_{\text{inj}}(\mc{A})$.
\begin{lem}\label{lem:exist_htpy_min}
Let $\mc{A}$ be a Grothendieck abelian category. Every complex $X^{\bullet}$ admits a quasi-isomorphism to a homotopically minimal $K$-injective complex with injective terms.
\end{lem}
\begin{proof}
The category $\mc{A}$ admits injective envelopes by \cite[Theorem~6.25]{Freyd}. By \cite{Stacks079P}, we can choose a quasi-isomorphism $X^{\bullet}\rightarrow I^{\bullet}$ where $I^{\bullet}$ is $K$-injective and each term $I^n$ is injective. By Proposition \ref{krause_dec}, there exists a decomposition $I^{\bullet}= J^{\bullet} \oplus N^{\bullet}$ such that $J^{\bullet}$ is homotopically minimal and $N^{\bullet}$ is nullhomotopic. The complex $J^{\bullet}$ is $K$-injective consisting of injective terms as it is a direct summand of $I^{\bullet}$ with the same properties. Since $N^{\bullet}$ is nullhomotopic, the composite $X^{\bullet}\rightarrow I^{\bullet}\twoheadrightarrow J^{\bullet}$ is the required quasi-isomorphism.
\end{proof}

\section{Admissibility of complexes}\label{sec:proof}

Let $G$ be a $p$-adic Lie group and $k$ be a field of characteristic $p$. Denote by $\mc{A}:= \text{Mod}_k(G)$ the abelian category of smooth $G$-representations in $k$-vector spaces and by $\text{Mod}_{\text{adm}}(G)$ its full subcategory of admissible representations. For a closed subgroup $H\subseteq G$, we denote by $\ind^G_H(-
)$ (resp. $\Ind^G_H(-)$) the compact induction (resp. full smooth induction) functors from $\text{Mod}_k(H)$ to $\text{Mod}_k(G)$. We denote by $D(G):=D(\mc{A})$ the unbounded derived category of $\mc{A}$.

We will use two standard facts: First, if $U\subseteq G$ is a
compact open pro-$p$ subgroup, then a smooth $G$-representation $V$
is admissible if and only if $V^U$ is finite-dimensional
(see \cite[Corollary 1.8 and the proof of Proposition 1.9]{Kohlhaase2017} or \cite[Lemma 2.22]{herzigaws}). Secondly, $\text{Mod}_{\text{adm}}(G)$ is a Serre subcategory of $\text{Mod}_k(G)$ (see, for example, \cite[Proposition 2.2.13]{Emerton2010I}), in particular, admissibility is preserved under subobjects, quotients and extensions. Also, $\mc{A}=\text{Mod}_k(G
)$ is a Grothendieck abelian category \cite[Lemma~1 (iv)]{Sch2015}. For $X^{\bullet}\in D(G)$ and a compact open subgroup $U\subseteq G$, we have the hypercohomology
\[
H^i(U, X^{\bullet}) := h^i(RH^0(U, X^{\bullet})).
\]
Note that the total right derived functor $RH^0(U,-):D(G)\rightarrow D(k)$ exists by \cite[Theorem~ I.5.1]{Hartshorne1966} using the triangulated subcategory $K_{\text{inj}}(\mc{A})\subseteq K(\mc{A})$.
We now prove our main theorem:
\begin{mainthm}\label{thm:main}
Let $U\subseteq G$ be a compact open pro-$p$ subgroup. Let $V^{\bullet}\in D(G)$ and $i\in \Z$. Then
\[
\dim_k H^i(U, V^{\bullet})< \infty \Rightarrow h^i(V^{\bullet})\; \text{is an admissible}\; G\text{-representation}.
\]
\end{mainthm}
\begin{proof}
Since restriction to $U$ is an exact functor which commutes with cohomology and admissibility can be checked by restricting to any compact open pro-$p$ subgroup of $G$, we can work with $\text{Res}^G_U(V^{\bullet})\in D(U)$, i.e., it suffices to assume $G=U$. Now, by Lemma \ref{lem:exist_htpy_min}, we can choose a homotopically minimal K-injective resolution with injective terms $V^{\bullet}\rightarrow J^{\bullet}$ in $\text{Mod}_k(U)$. Thus, the hypercohomology can be computed as
\[
H^i(U, V^{\bullet}) \cong h^i((J^{\bullet})^U).
\]
We now claim that the differentials vanish on the complex $(J^{\bullet})^U$. Indeed, fix $n\in \Z$, and let $v\in (J^n)^U$ with $v\neq 0$. Then $kv:=\{\lambda v: \lambda \in k\}$ is a subrepresentation of $J^n$. Since $Z^n(J^{\bullet})\hookrightarrow J^n$ is an injective envelope, essentiality implies that $kv \; \cap \ Z^n(J^{\bullet})\neq 0 \Rightarrow kv\subseteq Z^n(J^{\bullet})$. Thus, we have $d^n(v)=0$. Since $v \in (J^n)^U$ was arbitrary, we see that
\[
d^n((J^{n})^U)=0 \quad \forall \; n\in \Z.
\]
This means that
\begin{equation}\label{eq:main}
H^i(U, V^{\bullet}) \cong h^i((J^{\bullet})^U) = (J^i)^U. 
\end{equation}
Now, if $H^i(U, V^{\bullet})$ is finite dimensional, then $(J^i)^U$ is finite dimensional, i.e., $J^i$ is an admissible $U$-representation. Since admissibility is preserved under subquotients, the representation $h^i(J^{\bullet})$ is also an admissible $U$-representation. This completes the proof.
\end{proof}

To obtain equality of the two admissibility conditions, we recall Schneider--Sorensen's proof of the converse inclusion.
\begin{prop}[{\cite[Proposition 4.8]{SchSor}}]\label{prop:schsor}
$D_{\text{adm}}(G)\subseteq D(G)^a$.    
\end{prop}
\begin{proof}
Let $d=\dim_{\Q_p}(G)$. Let $U\subseteq G$ be a torsion-free compact open pro-$p$ subgroup. Then $U$ has $p$-cohomological dimension $d$. Thus, $RH^0(U, -)$ has cohomological amplitude $[0,d]$, in particular, 
\[
RH^0(U, D^{\geq a}(G))\subseteq D^{\geq a}(k),
\]
\[
RH^0(U, D^{\leq b}(G))\subseteq D^{\leq b+d}(k).
\]
Thus, the functor $RH^0(U,-)$ is way-out in both directions (in the sense of \cite[Chapter~I, \S7]{Hartshorne1966}). Moreover, if $V$ is an admissible representation, then $H^i(U,V)$ is finite dimensional for all $i\geq 0$ \cite[Lemma~4.7]{SchSor}. Since finite dimensional vector spaces form a Serre subcategory, the way-out lemma \cite[Proposition I.7.3(iii)]{Hartshorne1966} implies that $RH^0(U,-)$ sends $D_{\text{adm}}(G)$ to $D_{\text{fin}}(k)$, i.e., complexes with finite dimensional cohomology. This proves the required inclusion.
\end{proof}
\begin{proof}[Alternative proof via biduality]
For $X^{\bullet}\in D_{\text{adm}}(G)$, the canonical biduality morphism
\[
\eta_{X^{\bullet}}: X^{\bullet}\longrightarrow \mathbb{D}^2(X^{\bullet}), \quad \quad \mathbb{D}=R\underline{\Hom}(-,k)
\] is an isomorphism by \cite[Proposition~3.7]{SchSor}. This is equivalent to $X^{\bullet}\in D(G)^a$ by \cite[Theorem~4.5]{SchSor}, establishing the required inclusion.
\end{proof}

Combining this with
Theorem A then gives
\[
D_{\text{adm}}(G)=D(G)^a.
\]
\begin{rem}
The argument above also gives another proof for \cite[Proposition~4.11]{SchSor}. If $H^i(U, V^{\bullet})=0$, then \eqref{eq:main} gives $(J^i)^U=0$ which forces $J^i=0$ since $U$ is pro-$p$ \cite[Lemma 2.6]{herzigaws}. Thus, $h^i(V^{\bullet})=h^i(J^{\bullet})=0$. As
\[
\Hom_{D(G)}(\ind^G_U(k), V^{\bullet}[i])=H^i(U, V^{\bullet}),
\] this shows that $\mathbb{X}:=\ind^G_U(k)$ is a generator of $D(G)$. When $U$ is additionally torsion-free, Schneider shows that this generator is a compact object in $D(G)$ \cite[Lemma~4]{Sch2015}.
\end{rem}

\begin{rem}[The torsion-freeness assumption]\label{rem:torsion}
For a smooth representation $V \in \text{Mod}_k(G)$, finite-dimensionality of $V^U$ for
one compact open pro-$p$ subgroup $U$ implies admissibility, and hence
finite-dimensionality of $V^K$ for every compact open subgroup $K$;
no torsion-freeness assumption is needed. For unbounded complexes $X^{\bullet}\in D(G)$,
however, the condition
\[
             \dim_k H^i(U,X^\bullet)<\infty\qquad(i\in\Z)
\]
can depend on the choice of $U$ if subgroups with torsion are allowed.
Indeed, let $G=U$ be a compact pro-$p$ $p$-adic Lie group with torsion,
and choose a closed subgroup $C\cong\Z/p\Z$. Set $W=\Ind_C^U(k)$, i.e., $W$ consists of locally constant functions on $C\backslash U$. For every open subgroup $K\subseteq U$, the space
$W^K$ consists of functions on the finite set $C\backslash U/K\Rightarrow $ $W$ is an admissible $U$-representation. Shapiro's lemma gives
\[
                 H^i(U,W)\cong H^i(C,k)\cong k\qquad\forall\; i\ge0.
\]
Now take $X^\bullet=\bigoplus_{m\ge0}W[m]$ with zero differentials.
Then $X^\bullet\in D_{\adm}(U)$, but we have
an injection
\[
 \bigoplus_{m\ge0}k
 \cong\bigoplus_{m\ge0}H^m(U,W) \cong \bigoplus_{m\geq 0} H^0(U, W[m])
 \hookrightarrow H^0(U,X^\bullet).
\]
Thus $H^0(U,X^\bullet)$ is infinite-dimensional.
Choosing a torsion-free compact open pro-$p$ subgroup $U_0\subseteq U$,
on the other hand would imply that
$H^i(U_0,X^\bullet)$ is finite-dimensional for every $i$ by Proposition \ref{prop:schsor}.
The same complex therefore satisfies the finiteness condition for
$U_0$ but not for $U$.

The representation $W$ also shows that the total right derived functor $RH^0(U,-)$ (which exists by using $K$-injective resolutions) is not way-out
left as $W\in D^{\leq 0}(U)$ but $RH^0(U, W)\notin D^{\leq n}(k)$ for any $n\geq 0$. 
\end{rem}

\begin{rem}
Let $U$ be a $p$-torsion-free compact $p$-adic Lie group. 
The natural functor
\[
 F_U: D(\text{Mod}_{\adm}(U))\longrightarrow D_{\adm}(U)
\]
is an equivalence. This extends the bounded-below comparison of
\cite[Remark 4.13]{SchSor} to unbounded complexes under the additional $p$-torsion-freeness assumption. Indeed, first assume that $U$ is pro-$p$, hence torsion-free. The equality of the two
admissibility conditions and \eqref{eq:main} show that every object
of $D_{\adm}(U)$ can be represented by a $K$-injective complex with admissible
terms, proving essential surjectivity. Such a representative is also
$K$-injective in $\text{Mod}_{\adm}(U)$, since this is a full abelian
subcategory with exact inclusion. Computing morphisms into these
representatives in the homotopy category then proves full faithfulness
\cite[Lemma 13.31.2, Tag 070I]{Stacks070I}.

For a general $U$, choose an open pro-$p$ Sylow subgroup $P\subseteq U$. Then $P$ is torsion-free and $[U:P]$ is prime to $p$. Given $V^{\bullet}\in D_{\text{adm}}(U)$, the above argument gives a $K$-injective resolution
\[
\Res^U_P V^{\bullet}\xrightarrow[]{\sim} J^{\bullet}
\]
with admissible terms. Since $[U:P]$ is finite, full smooth induction preserves both admissibility and $K$-injectivity, i.e.,
\[
Q^{\bullet}:= \Ind^U_PJ^{\bullet}
\] is $K$-injective with admissible terms. Since $[U:P]$ is invertible in $k$, the unit of the adjunction $\Res^U_P \dashv \Ind^U_P$
\[
V^{\bullet}\rightarrow \Ind^U_P \Res^U_P V^{\bullet}
\] splits with left inverse given by 
the counit of the adjunction $\Ind^U_P \dashv \Res^U_P$
\[
\Ind^U_P \Res^U_P V^{\bullet}\rightarrow V^{\bullet}
\]
divided by $[U:P]$. Choose a homotopically minimal $K$-injective resolution $V^{\bullet}\rightarrow I^{\bullet}$ with injective terms. The splitting above then gives maps between $I^{\bullet}$ and $Q^{\bullet}$ in $D(U)$ whose composite is the identity on $I^{\bullet}$. Since both $I^{\bullet}$ and $Q^{\bullet}$ are $K$-injective, these maps can be represented by chain maps
\[
I^{\bullet}\xrightarrow[]{a}Q^{\bullet}\xrightarrow[]{b} I^{\bullet}, \quad ba\simeq \text{id}_{I^{\bullet}}.
\]
By homotopical minimality of $I^{\bullet}$, $ba$ is an isomorphism of chain complexes. Let $c:=(ba)^{-1}$ and $b':=cb$. Then $b'a = cba=\text{id}_{I^{\bullet}}$, i.e., the morphism of complexes 
\[
I^{\bullet}\xrightarrow[]{a}Q^{\bullet}
\] splits in every degree. Thus, $I^n$ is a direct summand of the admissible representation $Q^n$, and hence is admissible for every $n\in \Z$. This shows essential surjectivity of $F_U$. The full faithfulness follows from the same argument as for $U$ pro-$p$ above.

We cannot remove the $p$-torsion-freeness assumption. Indeed, take $U=\mathbb{Z}/p\mathbb{Z}$. Then admissible smooth $U$-representations are precisely finitely generated $k[U]$-modules and 
\[
k[U]\cong k[t]/(t^p)
\] 
is a quasi-Frobenius ring of infinite global dimension \cite[Section~3]{PositselskiSchnurer2021}. By \cite[Theorem~3.1]{PositselskiSchnurer2021}, the natural functor
\[
F_U: D(\text{Mod}_{\text{adm}}(U))\rightarrow D_{\text{adm}}(U)
\] is neither full nor essentially surjective.
\end{rem}

\section{An application to derived parabolic induction}
Let $F/\Q_p$ be a finite extension, and suppose $G=\mathbf{G}(F)$ is the group of $F$-points of a connected reductive group $\mathbf{G}$ over $F$. Fix a parabolic subgroup $P=MN$ and let $d=\text{dim}_{\Q_p}(N)$. First, our result immediately implies that derived parabolic induction $i^G_P:D(M)\rightarrow D(G)$ preserves global admissibility. Indeed, ordinary parabolic induction is exact and preserves admissible representations \cite[Proposition~4.2 and Corollary~4.7]{Vigneras2016}. Consequently, it's extension to derived categories satisfies
\[
i^G_P(D(M)^a)=i^G_P(D_{\text{adm}}(M)) \subseteq D_{\text{adm}}(G) = D(G)^a.
\]
This provides another proof of \cite[Proposition~4.2.6]{Heyer2023}.

Let 
\[
L(N,-):D(G)\rightarrow D(M)
\]
denote the left adjoint of derived parabolic induction as constructed by Heyer \cite{Heyer2023}. We can now give another proof of the following result:
\begin{prop}\cite[Proposition 4.2.12]{Heyer2023}
The functor $L(N,-)$ preserves global admissibility, i.e.,  
\[
L(N, D(G)^a)\subseteq D(M)^a.
\]
\end{prop}
\begin{proof}
The key representation-theoretic input is that Emerton's higher ordinary parts functor $H^j\text{Ord}_P(-)$ preserves admissible representations over arbitrary fields of characteristic $p$, as noted in \cite[Section~6]{HSS2019}. 

We also note that $L(N,-)$ has cohomological amplitude $[-d,0]$ by \cite[Proposition 4.1.5]{Heyer2023} and for $i\geq 0$, we have
\[
h^{-i}L(N,V)\cong H^{d-i}\text{Ord}_P(V)\otimes_k \delta_P
\] for an admissible smooth representation $V$ \cite[Remark 4.1.6]{Heyer2023}, where $\delta_P:P\rightarrow k^{\times}$ is the smooth character described in \cite[Remark 3.4.21]{Heyer2023}, viewed as a character on $M$. Admissibility of $H^{d-i}\text{Ord}_P(V)$ then implies that each $h^{-i}L(N,V)$ is an admissible $M$-representation, i.e., $L(N,V)\in D_{\text{adm}}(M)$. 

The finite cohomological amplitude also shows that $L(N,-)$ is way-out in both directions. Since admissible representations form a Serre subcategory, the way-out lemma \cite[Proposition I.7.3(iii)]{Hartshorne1966} gives
\[
L(N,D_{\text{adm}}(G))\subseteq D_{\text{adm}}(M).
\]
The equality $D_{\text{adm}}(H)=D(H)^a$ for $H=G,M$ now gives the desired preservation of global admissibility.

\end{proof}

\bibliographystyle{amsalpha}
\bibliography{references}

\end{document}